\documentclass[12pt]{article}

\usepackage{amssymb,amsfonts,latexsym}
\usepackage{pst-plot}

\newcommand{\ben}{\begin{enumerate}}
\newcommand{\een}{\end{enumerate}}
\newcommand{\barr}{\begin{array}}
\newcommand{\earr}{\end{array}}
\newcommand{\btab}{\begin{tabular}}
\newcommand{\etab}{\end{tabular}}
\newcommand{\beq}{\begin{equation}}
\newcommand{\eeq}{\end{equation}}
\newcommand{\bea}{\begin{eqnarray*}}
\newcommand{\eea}{\end{eqnarray*}}
\newcommand{\case}[4]
{\left\{\barr{ll}#1&\mbox{#2}\\#3&\mbox{#4}\earr\right.}
\newcommand{\iso}{\cong}
\newcommand{\Cong}{\equiv}
\newcommand{\Mod}{\mathop{\rm mod}\nolimits}

\newtheorem{theorem}{Theorem}[section]
\newtheorem{proposition}[theorem]{Proposition}
\newtheorem{lemma}[theorem]{Lemma}

\newcounter{claimcounter}
\newtheorem{claim}[claimcounter]{Claim}

\newcommand{\qed}{  \rule{1ex}{1ex}}
\newenvironment{proof}{\noindent {\bf Proof:}}{{\qed}}

\makeatletter
\newfont{\footsc}{cmcsc10 at 8truept}
\newfont{\footbf}{cmbx10 at 8truept}
\newfont{\footrm}{cmr10 at 10truept}
\title{Maximal Independent Sets In Graphs With At Most $r$ Cycles}
\author{Goh Chee Ying\\[-5pt]
\small Department of Mathematics\\[-5pt]
\small National University of Singapore\\[-5pt]
\small Singapore\\[-5pt]
\small \texttt{goh\_chee\_ying@moe.edu.sg}\\[6pt]
Koh Khee Meng\\[-5pt]
\small Department of Mathematics\\[-5pt]
\small National University of Singapore\\[-5pt]
\small Singapore\\[-5pt]
\small \texttt{matkohkm@nus.edu.sg}\\[6pt]
Bruce E. Sagan\\[-5pt]
\small Department of Mathematics\\[-5pt]
\small Michigan State University\\[-5pt]
\small East Lansing, MI\\[-5pt]
\small \texttt{sagan@math.msu.edu}\\[6pt]
Vincent R. Vatter\thanks{Partially 
supported by an award from DIMACS and
an NSF VIGRE grant to the Rutgers University Department of  
Mathematics.}\\[-5pt]
\small Department of Mathematics\\[-5pt]
\small Rutgers University\\[-5pt]
\small Piscataway, NJ\\[-5pt]
\small \texttt{vatter@math.rutgers.edu}}

\date{\today \\[6pt]
	\begin{flushleft}
	\end{flushleft}
           }

\begin{document}
\maketitle

\begin{abstract}
We find the maximum number of maximal independent sets in two families  
of graphs.  The first family consists of all graphs with $n$ vertices and at most $r$ cycles.
The second family is all graphs of the first family which are connected and satisfy $n\ge 3r$.
\end{abstract}

\section{Introduction}\label{mmi1-intro}

Let $G=(V,E)$ be a simple graph.  A subset $I\subseteq V$ is {\em  
independent\/} if there is no edge of $G$ between any two vertices of  
$I$.  Also, $I$ is {\em maximal\/} if it is not properly contained in  
any other independent set.  We let $m(G)$ be the number of maximal  
independent sets of $G$.

Around 1960, Erd\H{o}s and Moser asked for the maximum value of
$m(G)$ as $G$ runs over all graphs with $n$ vertices as well as for a
characterization of the graphs achieving this maximum.  (Actually, they
asked the dual question about cliques in such graphs.)  Shortly
thereafter Erd\H{o}s, and slightly later Moon and Moser~\cite{mm:cg},
answered both questions.  The extremal graphs turn out to have most of
their components isomorphic to the complete graph $K_3$.
Wilf~\cite{wil:nmi} raised the same questions for the family of
connected graphs.  Independently, F\H{u}redi~\cite{fur:nmi} determined
the maximum number for $n>50$, while Griggs, Grinstead, and
Guichard~\cite{ggg:nmi} found the maximum for all $n$ as well as
the extremal graphs.  Many of the blocks (maximal subgraphs
containing no cutvertex) of these graphs are also $K_3$'s.

Since these initial papers, there has been a string of articles about  
the maximum value of $m(G)$ as $G$ runs over various families of  
graphs.  In particular, graphs with a bounded number of cycles have  
received attention.  Wilf~\cite{wil:nmi} determined the maximum number  
of maximal independent sets possible in a tree, while  
Sagan~\cite{sag:nis} characterized the extremal trees.  These involve  
attaching copies of $K_2$ to the endpoints of a given path.  Later Jou  
and Chang~\cite{jc:mis} settled the problem for graphs and connected  
graphs with at most one cycle.  Here we consider the family of graphs
with $n$ vertices and at most $r$ cycles, and the family of connected graphs
with $n$ vertices and at most $r$ cycles where $n\ge 3r$.  The extremal graphs
are obtained by  
taking copies of $K_2$ and $K_3$ either as components (for all such  
graphs) or as blocks (for all such connected graphs).  We define the extremal
graphs and prove some lemmas about them in the next section.
Then Section~\ref{pmt} gives the proof of our main result,  
Theorem~\ref{main1}.

\section{Extremal graphs and lemmas}

For any two graphs $G$ and $H$, let $G\uplus H$ denote the disjoint
union of $G$ and $H$, and for any nonnegative integer $t$, let $tG$
stand for the disjoint union of $t$ copies of $G$.  We will need
the original result of Moon and Moser.  To state it, suppose $n\ge 2$ and
let
$$
G(n):=
\left\{\barr{ll}
\frac{n}{3}K_3&\mbox{if $n\Cong0\ (\Mod 3)$,}\\
2K_2\uplus\frac{n-4}{3}K_3&\mbox{if $n\Cong1\ (\Mod
3)$,}\rule{0pt}{20pt}\\
K_2\uplus\frac{n-2}{3}K_3&\mbox{if $n\Cong2\ (\Mod 3)$.}\rule{0pt}{20pt}
\earr\right.
$$
Also let
$$
\mbox{$G'(n):=K_4\uplus\frac{n-4}{3}K_3$ if $n\Cong1\ (\Mod 3)$.}
$$
Using the fact that $m(G\uplus H)=m(G)m(H)$
we see that
$$
g(n):=m(G(n))=
\left\{\barr{ll}
3^\frac{n}{3}&\mbox{if $n\Cong0\ (\Mod 3)$,}\\
4\cdot3^\frac{n-4}{3}&\mbox{if $n\Cong1\ (\Mod 3)$,}\rule{0pt}{20pt}\\
2\cdot3^\frac{n-2}{3}&\mbox{if $n\Cong2\ (\Mod 3)$.}\rule{0pt}{20pt}\\
\earr\right.
$$
Note also that if $n\Cong1\ (\Mod 3)$ then $m(G'(n))=m(G(n))$.

\begin{theorem}[Moon and Moser~\cite{mm:cg}]
\label{mm}
Let $G$ be a graph with $n\ge2$ vertices.  Then
$$
m(G)\le g(n)
$$
with equality if and only if $G\iso G(n)$ or, 
for $n\Cong1\ (\Mod 3)$,
$G\iso G'(n)$.\qed
\end{theorem}

Note that $G(n)$ has at most $\lfloor n/3\rfloor$ cycles.  Therefore the Moon-Moser Theorem gives the maximum number of maximal independent sets for the family of all graphs with $n$ vertices and at most $r$ cycles when $r\ge\lfloor n/3\rfloor$.  To complete the characterization, we need only handle the cases where $r<\lfloor n/3\rfloor$.  To make our proof cleaner, we will assume the stronger condition that $n\ge 3r-1$.

For any positive integers $n,r$ with $n\ge 3r-1$ we define
$$
G(n,r):=\case{rK_3\uplus\frac{n-3r}{2}K_2}
{if $n\Cong r\ (\Mod 2)$,}
{(r-1)K_3\uplus\frac{n-3r+3}{2}K_2}
{if $n\not\Cong r\ (\Mod 2)$.\rule{0pt}{20pt}}
$$

Note that if $n$ and $r$ have different parity then $G(n,r)\iso  
G(n,r-1)$.
This duplication is to facilitate the statement and proof of our main
result where $G(n,r)$ will be extremal among all graphs with $|V|=n$
and at most $r$ cycles.  Further, define
$$
g(n,r) := m(G(n,r)) =
\case{3^r \cdot 2^\frac{n-3r}{2}}
{if $n\Cong r\ (\Mod 2)$,}
{3^{r-1} \cdot 2^\frac{n-3r+3}{2}}
{if $n\not\Cong r\ (\Mod 2)$.\rule{0pt}{30pt}}
$$

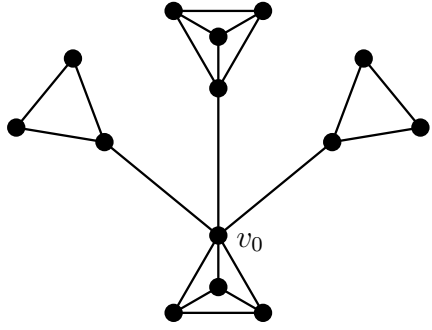
\begin{figure}
\begin{center}
\psset{xunit=0.012in, yunit=0.012in}
\psset{linewidth=1.0\psxunit}
\begin{pspicture}(0,0)(176.384,131.869)
\psline(151.631, 111.092)(176.384, 80.9303)     
\psline(88.1921, 120.606)(107.701, 131.869)     
\psline(68.6831, 131.869)(107.701, 131.869)     
\psline(68.6831, 131.869)(88.1921, 120.606)     
\psline(0, 80.9303)(24.7528, 111.092)   
\psline(137.887, 74.5744)(88.1921, 33.7906)     
\psline(137.887, 74.5744)(176.384, 80.9303)     
\psline(137.887, 74.5744)(151.631, 111.092)     
\psline(88.1921, 98.0785)(88.1921, 33.7906)     
\psline(88.1921, 98.0785)(107.701, 131.869)     
\psline(88.1921, 98.0785)(88.1921, 120.606)     
\psline(88.1921, 98.0785)(68.6831, 131.869)     
\psline(38.4969, 74.5744)(88.1921, 33.7906)     
\psline(38.4969, 74.5744)(24.7528, 111.092)     
\psline(38.4969, 74.5744)(0, 80.9303)   
\psline(68.6831, 0)(88.1921, 33.7906)   
\psline(107.701, 0)(88.1921, 33.7906)   
\psline(107.701, 0)(68.6831, 0) 
\psline(88.1921, 11.2635)(88.1921, 33.7906)     
\psline(88.1921, 11.2635)(68.6831, 0)   
\psline(88.1921, 11.2635)(107.701, 0)   
\pscircle*(88.1921, 33.7906){4\psxunit} 
\pscircle*(176.384, 80.9303){4\psxunit} 
\pscircle*(151.631, 111.092){4\psxunit} 
\pscircle*(107.701, 131.869){4\psxunit} 
\pscircle*(88.1921, 120.606){4\psxunit} 
\pscircle*(68.6831, 131.869){4\psxunit} 
\pscircle*(24.7528, 111.092){4\psxunit} 
\pscircle*(0, 80.9303){4\psxunit}       
\pscircle*(137.887, 74.5744){4\psxunit} 
\pscircle*(88.1921, 98.0785){4\psxunit} 
\pscircle*(38.4969, 74.5744){4\psxunit} 
\pscircle*(68.6831, 0){4\psxunit}       
\pscircle*(107.701, 0){4\psxunit}       
\pscircle*(88.1921, 11.2635){4\psxunit} 
\rput[c](102.1921, 30.7906){$v_0$}
\end{pspicture}
\end{center}
\caption{The graph $C(14)$}\label{C14}
\end{figure}

For the connected case, we will use the result of Griggs, Grinstead,  
and Guichard.
We obtain the extremal graphs as follows.
Let $G$ be a graph all of whose components are complete and let $K_m$
be a complete graph disjoint from $G$.  Construct the graph $K_m * G$ by
picking a vertex $v_0$ in $K_m$ and connecting it to a
single vertex in each component of $G$.  If $n\ge6$ then let
$$
C(n):=
\left\{\barr{ll}
K_3 * \frac{n-3}{3}K_3&\mbox{if $n\Cong0\ (\Mod 3)$,}\\
K_4 * \frac{n-4}{3}K_3&\mbox{if $n\Cong1\ (\Mod 3)$,}\rule{0pt}{20pt}\\
K_4 * \left(K_4\uplus\frac{n-8}{3}K_3\right)
		&\mbox{if $n\Cong2\ (\Mod 3)$.}\rule{0pt}{20pt}
\earr\right.
$$
The graph $C(14)$ is displayed in Figure~\ref{C14}.
Counting  maximal independent sets by whether they do or do not
contain $v_0$ gives
$$
c(n):=m(C(n))=
\left\{\barr{ll}
2 \cdot 3^\frac{n-3}{3} + 2^\frac{n-3}{3}&\mbox{if $n\Cong0\ (\Mod
3)$,}\\
3^\frac{n-1}{3} + 2^\frac{n-4}{3}
		&\mbox{if $n\Cong1\ (\Mod 3)$,}\rule{0pt}{20pt}\\
4 \cdot 3^\frac{n-5}{3} + 3 \cdot 2^\frac{n-8}{3}
		&\mbox{if $n\Cong2\ (\Mod 3)$.}\rule{0pt}{20pt}
\earr\right.
$$

\begin{theorem}[Griggs, Grinstead, and Guichard~\cite{ggg:nmi}]
\label{ggg}
Let $G$ be a connected graph with $n\ge6$
vertices.  Then
$$
m(G)\le c(n)
$$
with equality if and only if $G\iso C(n)$.\qed
\end{theorem}

In order to limit the number of cases in the proof of our main theorem we will only find the maximum of $m(G)$ for the family of all connected graphs when $n\ge 3r$.  Unlike the arbitrary graphs case, this result, together with the Griggs-Grinstead-Guichard Theorem, does not completely determine the maximum of $m(G)$ for all $n$ and $r$.  For example, when $n=10$ the extremal connected graph given by the Griggs-Grinstead-Guichard Theorem has 9 cycles, while our proof will only characterize extremal connected graphs with at most 3 cycles.  Although this gap between our main theorem and the Griggs-Grinstead-Guichard Theorem is relatively small ($C(n)$ contains $\lfloor n/3\rfloor$, $\lfloor n/3\rfloor + 6$, $\lfloor n/3\rfloor+12$ cycles when $n\Cong 0,1,2\ (\Mod 3)$ respectively), it takes considerable care to handle it.  This work is undertaken in \cite{mmi2}.

\begin{figure}
\begin{center}
\btab{ccc}
\psset{xunit=0.012in, yunit=0.012in}
\psset{linewidth=1.0\psxunit}
\begin{pspicture}(0,-33.7906)(184.776,133.7906)
\psline(145.841, 35.7164)(92.388, 0)    
\psline(116.99, 59.3943)(92.388, 0)     
\psline(92.388, 64.2879)(92.388, 0)     
\psline(67.786, 59.3943)(92.388, 0)     
\psline(38.9345, 35.7164)(92.388, 0)    
\psline(184.776, 38.2683)(145.841, 35.7164)     
\psline(163.099, 70.7107)(145.841, 35.7164)     
\psline(163.099, 70.7107)(184.776, 38.2683)     
\psline(130.656, 92.388)(116.99, 59.3943)       
\psline(92.388, 100)(92.388, 64.2879)   
\psline(54.1196, 92.388)(67.786, 59.3943)       
\psline(21.6773, 70.7107)(38.9345, 35.7164)     
\psline(0, 38.2683)(38.9345, 35.7164)   
\psline(0, 38.2683)(21.6773, 70.7107)   
\pscircle*(92.388, 0){4\psxunit}        
\pscircle*(145.841, 35.7164){4\psxunit} 
\pscircle*(116.99, 59.3943){4\psxunit}  
\pscircle*(92.388, 64.2879){4\psxunit}  
\pscircle*(67.786, 59.3943){4\psxunit}  
\pscircle*(38.9345, 35.7164){4\psxunit} 
\pscircle*(184.776, 38.2683){4\psxunit} 
\pscircle*(163.099, 70.7107){4\psxunit} 
\pscircle*(130.656, 92.388){4\psxunit}  
\pscircle*(92.388, 100){4\psxunit}      
\pscircle*(54.1196, 92.388){4\psxunit}  
\pscircle*(21.6773, 70.7107){4\psxunit} 
\pscircle*(0, 38.2683){4\psxunit}       
\rput[c](106.388, -3){$v_0$}
\end{pspicture}
&\rule{30pt}{0pt}&
\psset{xunit=0.012in, yunit=0.012in}
\psset{linewidth=1.0\psxunit}
\begin{pspicture}(0,0)(184.776,133.791)
\psline(145.841, 69.5071)(92.388, 33.7906)      
\psline(116.99, 93.1849)(92.388, 33.7906)       
\psline(92.388, 98.0785)(92.388, 33.7906)       
\psline(67.786, 93.1849)(92.388, 33.7906)       
\psline(38.9345, 69.5071)(92.388, 33.7906)      
\psline(184.776, 72.059)(145.841, 69.5071)      
\psline(163.099, 104.501)(145.841, 69.5071)     
\psline(163.099, 104.501)(184.776, 72.059)      
\psline(130.656, 126.179)(116.99, 93.1849)      
\psline(92.388, 133.791)(92.388, 98.0785)       
\psline(54.1196, 126.179)(67.786, 93.1849)      
\psline(21.6773, 104.501)(38.9345, 69.5071)     
\psline(0, 72.059)(38.9345, 69.5071)    
\psline(0, 72.059)(21.6773, 104.501)    
\psline(72.8789, 0)(92.388, 33.7906)    
\psline(111.897, 0)(92.388, 33.7906)    
\psline(111.897, 0)(72.8789, 0) 
\pscircle*(92.388, 33.7906){4\psxunit}  
\pscircle*(145.841, 69.5071){4\psxunit} 
\pscircle*(116.99, 93.1849){4\psxunit}  
\pscircle*(92.388, 98.0785){4\psxunit}  
\pscircle*(67.786, 93.1849){4\psxunit}  
\pscircle*(38.9345, 69.5071){4\psxunit} 
\pscircle*(184.776, 72.059){4\psxunit}  
\pscircle*(163.099, 104.501){4\psxunit} 
\pscircle*(130.656, 126.179){4\psxunit} 
\pscircle*(92.388, 133.791){4\psxunit}  
\pscircle*(54.1196, 126.179){4\psxunit} 
\pscircle*(21.6773, 104.501){4\psxunit} 
\pscircle*(0, 72.059){4\psxunit}        
\pscircle*(72.8789, 0){4\psxunit}       
\pscircle*(111.897, 0){4\psxunit}       
\rput[c](106.388, 30.7906){$v_0$}
\end{pspicture}
\\\\
C(13,2)
&&
C(15,3)
\etab
\end{center}
\caption{Examples of $C(n,r)$ for $n\ge 3r$}\label{Cex}
\end{figure}
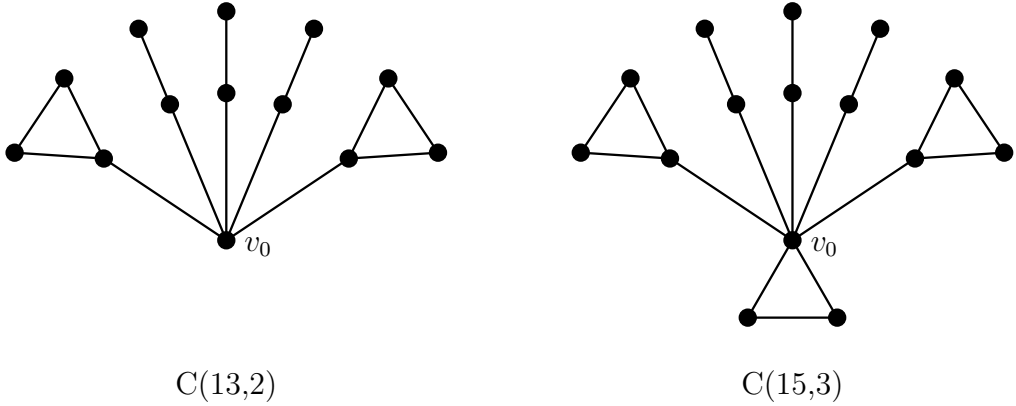

When $n\ge 3r$ we define
$$
C(n,r):=
\left\{\barr{ll}
K_3 * \left((r-1)K_3\uplus\frac{n-3r}{2}K_2\right)
		&\mbox{if $n\Cong r\ (\Mod 2)$,}\\
K_1 * \left(rK_3\uplus\frac{n-3r-1}{2}K_2\right)
		&\mbox{if $n\not\Cong r\ (\Mod 2)$.}\rule{0pt}{20pt}
\earr\right.
$$
The graphs $C(13,2)$ and $C(15,3)$ are shown
in Figure~\ref{Cex}.
As usual, we let
$$
c(n,r) := m(C(n,r)) =
\case{3^{r-1} \cdot 2^\frac{n-3r+2}{2} + 2^{r-1}}
{if $n\Cong r\ (\Mod 2)$,}
{3^r \cdot 2^\frac{n-3r-1}{2}}
{if $n\not\Cong r\ (\Mod 2)$.\rule{0pt}{30pt}}
$$

We also need the bounds for maximal independent sets in trees and  
forests, although we will not need the extremal graphs.  Define
$$
f(n):=2^{\lfloor\frac{n}{2}\rfloor}
$$
and
$$
t(n):=
\case{2^\frac{n-2}{2}+1}{if $n$ is even,}
{2^\frac{n-1}{2}}{if $n$ is odd.}
$$

Using our upcoming Proposition~\ref{m(G)}, 
it is easy to establish the following result.

\begin{theorem}
\label{forests}
If $G$ is a forest with $n\ge 1$ vertices then $m(G)\le f(n)$.\qed
\end{theorem}

Somewhat surprisingly, the tree analogue is significantly more  
difficult.

\begin{theorem}[Wilf~\cite{wil:nmi}]
\label{trees}
If $G$ is a tree with $n$ vertices then $m(G)\le t(n)$.\qed
\end{theorem}

For the extremal trees the reader is referred to Sagan~\cite{sag:nis}.

Next, we have a list of inequalities that will be useful in the proof  
of our main theorem.  It will be convenient to let $g(n,0)=f(n)$ and $c(n,0)=t(n)$.
\begin{lemma}\label{inequalities}
We have the following monotonicity results.
\ben
\item[(1)]\label{gnm}
If $r\ge1$ and $n>m\ge3r-1$ then
	$$g(n,r)>g(m,r).$$
\item[(2)]\label{cnm}
If $r\ge1$ and $n>m\ge3r$ then
	$$c(n,r)>c(m,r).$$
\item[(3)]\label{grq}
If $r>q\ge0$ and $n\ge3r-1$ then
	$$g(n,r)\ge g(n,q)$$
with equality if and only if $n$ and $r$ have different parity and
$q=r-1$.
\item[(4)]\label{crq}
If $r>q\ge0$ and $n\ge3r$ then
	$$c(n,r)\ge c(n,q)$$
with equality if and only if $(n,r,q)=(4,1,0)$ or $(7,2,1)$.
\een
\end{lemma}
\begin{proof}
The proofs of all of these results are similar, so we will content
ourselves with a demonstration of (4).  It  suffices
to consider the case when $q=r-1$.  Suppose that $r\ge 2$ since the $r=1$ case
is similar.  If $n$ and $r$ have the same
parity, then we wish to show
$$
3^{r-1} \cdot 2^\frac{n-3r+2}{2} + 2^{r-1}>
3^{r-1} \cdot 2^\frac{n-3r+2}{2}
$$
which is clear.  If $n$ and $r$ have different parity, then
$n\ge3r$ forces $n\ge3r+1$.  We want
$$
3^r \cdot 2^\frac{n-3r-1}{2} \geq
3^{r-2} \cdot 2^\frac{n-3r+5}{2} +2^{r-2}.
$$
Combining the terms with powers of 3, we have the equivalent
inequality
$$
3^{r-2} \cdot 2^\frac{n-3r-1}{2} \geq2^{r-2}.
$$
The bounds on $n$ and $r$ show that this is true, with equality exactly
when both sides equal 1.
\end{proof}

We also need two results about $m(G)$ for general graphs $G$.  In what  
follows, if $v\in V$ then the {\em open\/} and {\em closed  
neighborhoods of $v$} are
$$
N(v)=\{u\in V\ |\ uv\in E\}
$$
and
$$
N[v] = \{v\}\cup N(v),
$$
respectively.
We also call a block an
{\em endblock of $G$\/} if it has at most one cutvertex in the
graph as a whole.  We first verify that certain types of endblocks  
exist.

\begin{proposition}
\label{endB}
Every graph has an endblock that intersects at most one non-endblock.
\end{proposition}
\begin{proof}
The {\em block-cutvertex graph of $G$}, $G'$, is the
graph with a vertex $v_B$ for each block $B$ of $G$, a vertex $v_x$
for each cutvertex $x$ of $G$, and edges of the form $v_Bv_x$ whenever
$x\in V(B)$.  It is well known that $G'$ is a forest.  Now consider a
longest path $P$ in $G'$.  The  final vertex of $P$ corresponds to a
block $B$ of $G$ with the desired property.
\end{proof}

Any block with at least 3 vertices is 2-connected, i.e., one
must remove at least 2 vertices to disconnect or trivialize the graph.
Such graphs are exactly those which can be obtained from a cycle
by adding a sequence of {\it ears\/}.  This fact is originally due to
Whitney~\cite{w:ear}, and can also be found in Diestel~\cite[Proposition 3.1.2]{die:gt}
and West~\cite[Theorem 4.2.8]{w:gt}.
\begin{theorem}[Ear Decomposition Theorem]
A graph $B$ is 2-connected if and only if there is a sequence
$$
B_0, B_1,\ldots, B_l=B
$$
such that $B_0$ is a cycle and $B_{i+1}$ is obtained  by
taking a nontrivial path and identifying its two endpoints with two
distinct
vertices of $B_i$. \quad\qed
\end{theorem}

\begin{proposition}\label{m(G)}
The invariant $m(G)$ satisfies the following.
\ben
\item[(1)]\label{G-v}
If $v\in V$ then
$$
m(G)\le m(G-v)+m(G-N[v]).
$$
\item[(2)]\label{endblock}
If $G$ has an endblock $B$ that is isomorphic 
to a complete graph
then
$$
m(G)=\sum_{v\in V(B)} m(G-N[v]).
$$
\een
In fact, the same equality holds for any complete subgraph $B$
having at least one vertex that is adjacent in $G$ only to other
vertices of $B$.
\end{proposition}
\begin{proof}
For any $v\in V$ there is a bijection between the maximal
independent sets $I$ of $G$ that contain $v$ and the maximal
independent sets of $G-N[v]$, given by $I\mapsto I-v$.  Also, the
identity
map gives an injection from those $I$ that do not contain $v$ into
the maximal independent sets of $G-v$.  This proves (1).
For (2), merely use the previous bijection and the fact
that, under either hypothesis, any maximal independent set of $G$ must
contain exactly one of the vertices of $B$.
\end{proof}

We will refer to the formulas in parts (1) and (2) of this proposition as the
{\it $m$-bound\/} and {\it $m$-recursion\/}, respectively.

\section{Proof of the main theorem}\label{pmt}

We are now in a position to state and prove our main result.
The path and cycle on $n$ vertices will be denoted by $P_n$ and $C_n$,
respectively.  Also, let $E$ denote the graph pictured in
Figure~\ref{exception}.

\begin{theorem}\label{main1}
Let $G$ be a graph with $n$ vertices and at most $r$ cycles where
$r\ge1$.
\ben
\item[(I)]
If $n\ge3r-1$ then for all such graphs we have
$$m(G)\le g(n,r)$$
with equality if and only if $G\iso G(n,r)$.
\item[(II)]
If $n\ge3r$ then for all such graphs that are connected we have
$$m(G)\le c(n,r)$$
with equality if and only if $G\iso C(n,r)$ or if $G$ is one of the
exceptional cases listed in the following table.
$$
\barr{c|c|c}
n	&r	&\mbox{possible $G\not\iso C(n,r)$}\\
\hline
4	&1	& P_4\\
5	&1	& C_5\\
7	&2	& C(7,1), E
\earr
$$
\een
\end{theorem}
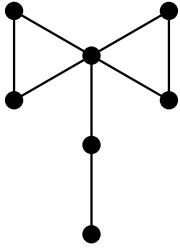
\begin{figure}
\begin{center}
\psset{xunit=0.012in, yunit=0.012in}
\psset{linewidth=1.0\psxunit}
\begin{pspicture}(0,0)(67.5813,97.5452)
\psline(67.5813, 97.5452)(33.7906, 78.0361)     
\psline(67.5813, 58.5271)(33.7906, 78.0361)     
\psline(67.5813, 58.5271)(67.5813, 97.5452)     
\psline(33.7906, 39.0181)(33.7906, 78.0361)     
\psline(0, 58.5271)(33.7906, 78.0361)   
\psline(0, 97.5452)(33.7906, 78.0361)   
\psline(0, 97.5452)(0, 58.5271) 
\psline(33.7906, 0)(33.7906, 39.0181)   
\pscircle*(33.7906, 78.0361){4\psxunit} 
\pscircle*(67.5813, 97.5452){4\psxunit} 
\pscircle*(67.5813, 58.5271){4\psxunit} 
\pscircle*(33.7906, 39.0181){4\psxunit} 
\pscircle*(0, 58.5271){4\psxunit}       
\pscircle*(0, 97.5452){4\psxunit}       
\pscircle*(33.7906, 0){4\psxunit}       
\end{pspicture}
\end{center}
\caption{The exceptional graph $E$}\label{exception}
\end{figure}
\begin{proof}  The proof will be by double induction on $n$ and $r$.   
The base
case $r=1$ has been done by Jou and Chang~\cite{jc:mis}, so we assume  
from now on
that $r\ge2$.  We will also assume that $n\ge 8$, as the smaller cases  
have been checked by computer.

We first show that graphs with a certain cycle structure can't be
extremal by proving the following pair of claims.  Here we assume that
$G$ has $n$ vertices and at most $r$ cycles.

\begin{claim}
If $G$ is a graph with two or more intersecting cycles and $n\ge3r-1$  
then $m(G)<g(n,r)$.
\end{claim}

\begin{claim}
If $G$ is a connected graph with an endblock $B$ containing two or more
cycles and $n\ge3r$ then $m(G)< c(n,r)$.
\end{claim}

To prove Claim 1 suppose to the contrary that $v$ is a vertex where two
cycles intersect, so $G-v$ has $n-1$ vertices and at most $r-2$ cycles.  
Furthermore, among all such vertices we can choose $v$ with
$\deg v\ge 3$.  It follows that $G-N[v]$ has at most $n-4$ vertices and
at most $r-2$ cycles.  If $r=2$ then using Theorem~\ref{forests} and the
$m$-bound gives
\bea
m(G)&\le&f(n-1)+f(n-4)\\[10pt]
&=&\case{2^\frac{n-2}{2}+2^\frac{n-4}{2}}{if $n$ is even,}
	{2^\frac{n-1}{2}+2^\frac{n-5}{2}}{if $n$ is odd}\\[10pt]
&=&\case{3\cdot 2^\frac{n-4}{2}}{if $n$ is even,}
	{5\cdot 2^\frac{n-5}{2}}{if $n$ is odd}\\[10pt]
&<&g(n,2).
\eea
If $r\ge 3$ then we use the induction hypothesis of the theorem,  
Lemma~\ref{inequalities}~(1) and~(3), and the $m$-bound to
get
\bea
m(G)&\le&g(n-1,r-2)+g(n-4,r-2)\\[10pt]
        &=&\case{3^{r-3} \cdot 2^\frac{n-3r+8}{2}
	+3^{r-2} \cdot 2^\frac{n-3r+2}{2} }
	{if $n\Cong r\ (\Mod 2)$,}
	{3^{r-2} \cdot 2^\frac{n-3r+5}{2}
	+3^{r-3} \cdot 2^\frac{n-3r+5}{2} }
	{if $n\not\Cong r\ (\Mod 2)$,\rule{0pt}{30pt}}\\[10pt]
        &=&\case{22\cdot3^{r-3} \cdot 2^\frac{n-3r}{2} }
	{if $n\Cong r\ (\Mod 2)$,}
	{8\cdot3^{r-3} \cdot 2^\frac{n-3r+3}{2} }
	{if $n\not\Cong r\ (\Mod 2)$,\rule{0pt}{30pt}}\\[10pt]
        &<&g(n,r).
\eea

To prove Claim 2 we first observe that by the Ear Decomposition Theorem, $B$  
contains two cycles $C$
and $C'$ such that $C\cap C'$ is a path with at least 2 vertices.
Since $B$ is an endblock, one of the endpoints of the path $C\cap C'$
is not a cutvertex in $G$.  Label this endpoint $v$.  Note that by  
construction $\deg v\ge3$ and  $v$ is on at least 3
cycles of $G$, namely $C$, $C'$, and the cycle in 
$C\cup C'$ obtained by
not taking any edge of the path $C\cap C'$.  Proceeding as in the
proof of Claim 1 and noting that $G-v$ is connected by our choice of  
$v$, we have
\begin{eqnarray*}
m(G)
&\le&
\case{t(n-1)+f(n-4)}{if $r=3$,}{c(n-1,r-3)+g(n-4,r-3)}{if $r>3$}\\
&<&
c(n,r).
\end{eqnarray*}

We now return to the proof of the theorem,  first tackling the case  
where $G$ varies over all graphs with $n$ vertices and at most $r$  
cycles.  For the base cases of $n=3r-1$ or $3r$, we have $g(n,r)=g(n)$  
and $G(n,r)=G(n)$ so we are done by the Moon-Moser Theorem.

Suppose that $n\ge3r+1$.  From Claim 1 we can assume that the cycles of
$G$ are disjoint.  It follows that the blocks of $G$ must all be cycles
or copies of $K_2$.  Let $B$ be an endblock of $G$.  We have three  
cases depending on whether $B\iso K_2$, $K_3$, or $C_p$ for $p\ge4$.

If $B\iso K_2$ then let $V(B)=\{v,w\}$ where $w$ is the cutvertex of  
$B$
in $G$, if $B$ has one.  Then $G-N[v]$ has $n-2$
vertices and at most $r$
cycles while $G-N[w]$ has at most $n-2$ vertices and at most $r$
cycles.  By induction, Lemma~\ref{inequalities}~(1) and~(3), and
the $m$-recursion we have
$$
m(G)\le 2g(n-2,r)=g(n,r)
$$
with equality if and only if $G-N[v]=G-N[w]\iso G(n-2,r)$.  It follows
that $B$ is actually a component of $G$ isomorphic to $K_2$ and so
$G\iso G(n,r)$.

The case $B\iso K_3$ is similar.
Proceeding as before, one obtains
$$
m(G)\le 3g(n-3,r-1)=g(n,r),
$$
and equality is equivalent to $G\iso B\uplus G(n-3,r-1)\iso G(n,r)$.

To finish off the induction step, consider $B\iso C_p$, $p\ge4$.  Then
there exist $v,w,x\in V(B)$ all of degree 2 such that $vw, vx\in E(B)$.
So $G-v$ has $n-1$ vertices and at most $r-1$ cycles.  Furthermore,
$G-v\not\iso G(n-1,r-1)$ since $G-v$ contains two vertices, $w$ and  
$x$, both of degree 1 and
in the same component but not adjacent.  Also, $G-N[v]$ has $n-3$
vertices and at most $r-1$ cycles.  Using computations similar to
those in Claim 1,
$$
m(G)<g(n-1,r-1)+g(n-3,r-1)=g(n,r),
$$
so these graphs cannot be extremal.

It remains to consider the connected case.  It will be convenient to  
leave the base cases of $n=3r$ or $3r+1$ until last, so assume that  
$n\ge 3r+2$.  Among all the endblocks of the form guaranteed by  
Proposition~\ref{endB}, let $B$ be one with the largest number of  
vertices.  Claim~2 shows that $B$ is either $K_2$ or a cycle.   
Furthermore, the $r=1$ base case 
shows that cycles with more  
than $5$ vertices are not extremal, so $B$ must contain a cutvertex  
$x$.  Again, there are three cases depending on the nature of $B$.

If $B\iso K_2$ then let $V(B)=\{x,v\}$  so that $\deg v=1$ and $\deg  
x\ge2$.  By the choice of $B$, $G-N[v]$ is the union of some number of  
$K_1$'s and a connected graph with at most $n-2$ vertices and at most  
$r$ cycles.  Also, $G-N[x]$ has at most $n-3$ vertices and at most $r$  
cycles, so
$$
m(G)\le c(n-2,r)+g(n-3,r)=c(n,r)
$$
with equality if and only if both $G-N[v]$ and $G-N[x]$ are extremal.   
Except for the case where $n=9$ and $r=2$, this implies that  
$G-N[v]\iso C(n-2,r)$ and $G-N[x]\iso G(n-3,r)$, which is equivalent to
$G\iso C(n,r)$.  In the case where $n=9$ and $r=2$ we still must have  
$G-N[x]\cong G(6,2)\cong 2K_3$, but now there are three possibilities  
for $G-N[v]$: $C(7,2)$, $C(7,1)$, or $E$.  However, since $G-N[x]$ is a
subgraph of $G-N[v]$ we must have $G-N[v]\iso C(7,2)$, which shows that
$G\iso C(9,2)$, as desired.

Next consider $B\iso K_3$ and let $V(B)=\{x,v,w\}$ where $x$ is the  
cutvertex.  Let $i$ be the  
number of $K_3$ endblocks other than
$B$ containing $x$.  First we note that $x$ is adjacent to some vertex
$y$ not in a $K_3$ endblock as otherwise $n<3r$.  It follows
from our choice of $B$ that $G-N[v]=G-N[w]$ has some number of
$K_1$ components, $i$ components isomorphic to $K_2$, and at most one
other component, say $H$, with at most $n-2i-3$ vertices and at most  
$r-i-1$
cycles.  Furthermore, because $x$ is adjacent to $y$, the graph
$G-N[x]$ has at most $n-2i-4$ vertices and at most $r-i-1$ cycles.   
This gives us the upper bound
$$
m(G)\le 2^{i+1}c(n-2i-3,r-i-1)+g(n-2i-4,r-i-1).
$$
As the right-hand side of this inequality is strictly decreasing for $i$ of a given parity, it suffices to consider the cases where $i$ is $0$ or $1$.  When $i=1$, we have
$$
m(G)\le 4c(n-5,r-2)+g(n-6,r-2)<c(n,r).
$$
When $i=0$ we have
$$
m(G)\le 2c(n-3,r-1)+g(n-4,r-1)=c(n,r).
$$
Using the same argument as in the case $B\cong K_2$ and
$(n,r)=(9,2)$, one can show that this inequality is strict when $c(n-3,r-1)$
could be achieved by one of the exceptional graphs.  For other $n,r$ we
get equality if and only if $G\cong C(n,r)$.

The last case is where $B\iso C_p$ where $p\ge4$.  Label the vertices of
$B$ as $x,u,v,w,\ldots$ so that they read one of the possible directions
along the cycle, where $x$ is the cutvertex.  So 
$\deg u=\deg v=\deg w =2$, $\deg x\ge3$, and $G-v$ is connected with $n-1$ vertices and at  
most $r-1$ cycles.  Furthermore, $G-
v\not\iso C(n-1,r-1)$ because it contains a vertex of degree 1
(namely $u$) adjacent to a vertex of degree at least 3 (namely $x$).
Also, the graph $G-N[v]$ is connected with $n-3$ vertices and at most
$r-1$ cycles.  These conditions give us
$$
m(G)<c(n-1,r-1)+c(n-3,r-1)=c(n,r),
$$
and since this inequality is strict, such $G$ are not extremal.

We are left with the base cases.  When $n=3r$, $c(n,r)=c(n)$ and
$C(n,r)\iso C(n)$ so we are done by the Griggs-Grinstead-Guichard  
Theorem.
If $n=3r+1$ then we can proceed as in the induction step except where  
$B\iso K_2$ since then $c(n-2,r)$ and $g(n-3,r)$ have arguments outside
of the permissible range.  However since we have assumed $n\ge 8$,  
Theorems~\ref{mm} and~\ref{ggg} apply to give
\bea
m(G) 	&\le& c(n-2)+g(n-3)\\
		&=& c(3r-1)+g(3r-2)\\
&=& 4 \cdot 3^{r-2} + 3 \cdot 2^{r-3} +  4 \cdot 3^{r-2}\\
&\le&c(3r+1,r)\\
&=&c(n,r).
\eea
The latter inequality is strict for $r\neq 3$; when $r=3$,
the former
inequality is strict because $C(8)=K_4*K_4$ has more than $3$ cycles.
Therefore these graphs cannot be extremal, finishing the proof
of the theorem.
\end{proof}

As was mentioned in Section~\ref{mmi1-intro}, Theorem~\ref{main1} and the Moon-Moser Theorem combine to completely settle the maximal independent set question for the
family of arbitrary graphs with $n$ vertices and at most $r$ cycles
for all $n$ and $r$, but this does
not occur in the connected case.  For these graphs,
Theorem~\ref{main1} handles the cases where $n$ is large
relative to $r$ and the Griggs-Grinstead-Guichard Theorem
handles the cases where $n$ is small relative to $r$, but there
is a gap between where these two results apply when
$n\not\Cong 0\ (\Mod 3)$.  This gap is handled in \cite{mmi2}.  That paper also
answers related questions for maximum independent sets.


\begin{thebibliography}{99}

\bibitem{die:gt} R. Diestel, ``Graph Theory,'' Graduate Texts in
Math., Vol.\ 173, Springer-Verlag, New York, NY, 1997.

\bibitem{fur:nmi} Z. F\H{u}redi, The number of maximal independent sets
in connected graphs, {\it J. Graph Theory\/} {\bf 11} (1987),
463--470.

\bibitem{ggg:nmi}  J. R. Griggs, C. M. Grinstead, and D. R. Guichard,
The number of maximal independent sets in a connected graph,
{\it Discrete Math.\/} {\bf 68} (1988), 211--220.

\bibitem{jc:mis} M. J. Jou and G. J. Chang, Maximal independent sets
in graphs with at most one cycle, {\it Discrete Appl.\ Math.\/}
{\bf 79} (1997), 67--73.

\bibitem{mm:cg} J. W. Moon and L. Moser, On cliques in graphs,
{\it Isr.\ J. Math.\/} {\bf 3} (1965), 23--28.

\bibitem{sag:nis} B. E. Sagan, A note on independent sets in trees,
{\it SIAM J. on Discrete Math.\/} {\bf 1} (1988), 105--108.

\bibitem{mmi2} B. E. Sagan and V. R. Vatter, Maximal and maximum independent sets in graphs with at most $r$ cycles, {\it J. Graph Theory\/}, to appear.

\bibitem{w:gt} D. West, ``Introduction to Graph Theory,'' 2nd edition, Prentice Hall, Upper Saddle River, NJ, 2001.

\bibitem{w:ear} H. Whitney, Congruent graphs and the connectivity of graphs, {\it Amer. J. Math.\/} {\bf 54} (1932), 150--168.

\bibitem{wil:nmi} H. S. Wilf, The number of maximal independent sets
in a tree, {\it SIAM J. Alg.\ Discrete Methods\/} {\bf 7} (1986),
125--130.

\end{thebibliography}
\end{document}